## RESEARCH ARTICLE

## THE TRANSFORM OF A LINE OF DESARGUES AFFINE PLANE IN AN ADDITIVE GROUP OF ITS POINTS

### [1],*Orgest Zaka and [2]Prof.Dr. Kristaq Filipi

[1]Department of Mathematics, Faculty of Technical Science, University of Vlora "Ismail QEMALI", Vlora, Albania
[2]Department of Mathematics, Faculty of Mathematical and Physical Engineering, Polytechnic University of Tirana, Tirana, Albania

---

**ABSTRACT**

In this paper we present a set transformation of points in a line of the Desargues affine plane in a additive group. For this, the first stop on the meaning of the Desargues affine plane, formulating first axiom of his that show proposition D1. Afterwards we show that little Pappus theorem, which we use in the construction of group proofs in additions of points on a line on desargues plane, also applies in the Desargues affine plane.

---

## INTRODUCTION

### Desargues affine plane

**Definition 1.1.** (Francis Borceux, 2014; Orgest ZAKA, Kristaq FILIPI 2016) *Affine plane called the incidence structure* $\mathcal{A} = (\mathcal{P}, \mathcal{L}, I)$ *that satisfies the following axioms:*

**A1:** *For every two different points $P$ and $Q \in \mathcal{P}$, there exists exactly one line $\lambda \in \mathcal{L}$ incident with that points.*
The line $\ell$, determined from the point $P$ and $Q$ will denoted $PQ$.
**A2:** *For a point $P \in \mathcal{P}$, and an line $\lambda \in \mathcal{L}$ such such $(P, \ell) \in \mathcal{I}$, there exists one and only one line line $r \in \mathcal{L}$, incident with point $P$ that such that $\ell \cap r = \varnothing$.*
**A3:** *In $\mathcal{A}$ there are three non-incident points with a line..*

The fact $(P, \ell) \in I$ (equivalent to P I $\ell$) we mark $P \in \ell$ and read *point $P$ is incident with a line $\ell$ or a line $\ell$ passes by points $P$ (contains point $P$)*. Whereas a straight line of the affine plane we consider as sets of points of affine plane with her incidents. From axioms **A1** implicates that tow different lines of $\mathcal{L}$ many have an common point, in other words *tow different lines of $\mathcal{L}$ or no have common point or have only one common point.*

---

*\*Corresponding author: Orgest Zaka,*
Department of Mathematics, Faculty of Technical Science, University of Vlora "Ismail QEMALI", Vlora, Albania.

**Definition 1.2.** *Two lines $\ell, m \in \mathcal{L}$ that matching or do not have common point of called parallel and in this case write $\ell \| m$, and when they have only one common point say that they expected.*

For single line $r \in \mathcal{L}$, which passes by a point $P \in \mathcal{P}$ and it is parallel with line $AB$, that does not pass the point $P$, we will use the notation $\ell_{AB}^{P}$.

PROPOSITION 1.1. (SADIKI, 2015) *Parallelism relation $\| = \{(r, s) \in \mathcal{L}^2 | r \| s\}$ on $\mathcal{L}$ is an equivalence relation in $\mathcal{L}$.*

**Definition 1.3.** *Three different points $P, Q, R \in \mathcal{P}$ wi called collineary, if there are incidents with the same line.*

**Definition 1.4.** *The set of three different non-collineary points $A$, $B$, $C$ together with the line $AB$, $BC$, $CA$ called three-vertex and marked $ABC$. Points $A$, $B$, $C$ called vertices, while the line $AB$, $BC$, $CA$ called side of three-vertex $ABC$.*

In affine Euclidian plane is true this

PROPOSITION D1. (Axiom I of Desargues) *If $AA_1$, $BB_1$, $CC_1$ are the three different parallel line (Fig. 1), then*

$$\left. \begin{array}{c} AB \| A_1 B_1 \\ BC \| B_1 C_1 \end{array} \right] \Rightarrow AC \| A_1 C_1.$$

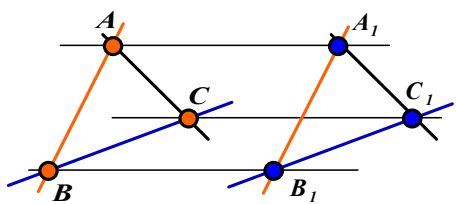

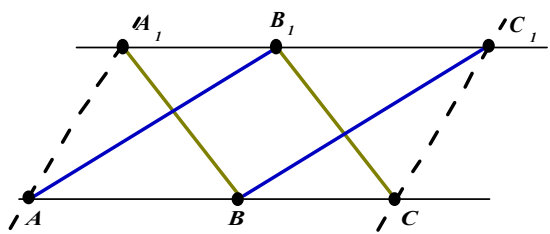

**Fig. 1.Desargues configuration**          **Fig. 2.Pappus configuration**

There are affine plans where Propositioni D1 not valid. Such a is the Moulton plane (9).

**Definition 1.5.** (Francis Borceux, 2014; SADIKI, 2015; COXETER, 1969) *An affine plane complete with desargues axiom D1, we shall call Desargues affine plane.*

Let's be now $A$, $B$, $C$ three different points of a line and $A_1$, $B_1$, $C_1$ three different points of one another straight-parallel to the first (Fig.2). If $AB_1 \| BC_1$ and $A_1 B \| B_1 C$ can contend that also $AA_1 \| CC_1$ ? Otherwise, we add the problem if it is true that

PROPOSITION 1.2 (FRANZ ROTHE, 2010; ROBÌN HARTSHORNE, 2000) (Little Pappus Theorem). *Let us be $A, B, C$ and $A_1, B_1, C_1$ two triple point located in two parallel lines (Fig. 2). If $AB_1 \| BC_1$ and $BA_1 \| CB_1$, then we have to $AA_1 \| CC_1$.*

The answer is that

THEOREM 1.1 (FRANZ ROTHE, 2010) (*the little Hessenberg Theorem*). **In the Desargues plane is tru** *Propositions 1.2, to wit is worth the Little Pappus theorem.*

*Proof.* Let us have two triplets of points $A$, $B$, $C$ and $A_1, B_1, C_1$ in two parallel lines such that $AB_1 \| BC_1$ and $BA_1 \| CB_1$ (Fig. 3).

We build a line $\ell_{AB_1}^{C}$ (a line that passes through points $C$ and it is parallel to the line $AB_1$), and line $\ell_{BA_1}^{A}$ (a line that passes through points $A$ and it is parallel to the line $BA_1$). We mark $D = \ell_{AB_1}^{C} \cap \ell_{BA_1}^{A}$. Also construct the line $DB$ (Fig.3).

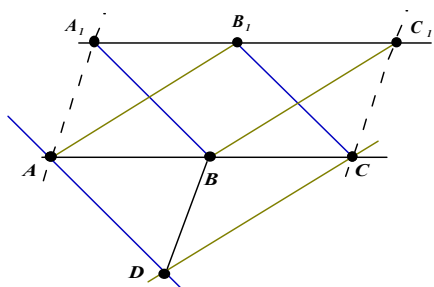

**Fig. 3. The proof configuration.**

By condition  of Propositions 1.2 we have $AB_1 \| BC_1$ and $BA_1 \| CB_1$. This imply parallelism of straight lines $AB_1$, $BC_1$, $\ell^C_{AB_1}$, and parallelism of straight lines $BA_1$, $CB_1$, $\ell^A_{BA_1}$. In these conditions, three-vertex $CC_1B_1$ and $DBA$ they have vertices in different parallel  line $AB_1$, $BC_1$, $\ell^C_{AB_1}$ and their sides satisfy the condition $AB \| B_1C_1$ and $AD \| B_1C$. Hence, according to axiom D1, we have also $CC_1 \| DB$.

Also, three-vertex $AA_1B_1$ and $DBC$ they have vertices in different parallel  line $BA_1$, $CB_1$, $\ell^A_{BA_1}$ and their sides satisfy the condition $BC \| A_1B_1$ and $DC \| AB_1$. Hence, according to axiom D1, we have also $DB \| AA_1$. By comparing the two conclusions of the implementation of axiom D1, according to Propositions 1.1, we conclude $AA_1 \| CC_1$.

## 2. Equipment of sets of points to a straight lines of the desargues affine plane with binary additive operations

In an Desargues affine plane $\mathcal{D} = (\mathcal{P}, \quad, \mathrm{I})$ we fix two different points  $O, I \in \mathcal{P}$ , which, according to axiom A1, determine a line $OI \in$    . Let us be $A$ and $B$  two whatever points of a  line $OI$. Choosing  in plane D  a point $B_1$ non-incidents with $OI$: $B_1 \notin OI$. Construct  line $\ell^{B_1}_{OI}$, which is only by axiom A2. Then construct  line $\ell^A_{OB_1}$, which also is the only according to axiom A2. Marking their intersection $P_1 = \ell^{B_1}_{OI} \cap \ell^A_{OB_1}$. Finally construct  line $\ell^{P_1}_{BB_1}$. For as much as $BB_1$ expects $OI$ in point $B$, then this straight line, parallel with $BB_1$, expects  line $OI$ in a single point $C$ (Fig. 4).

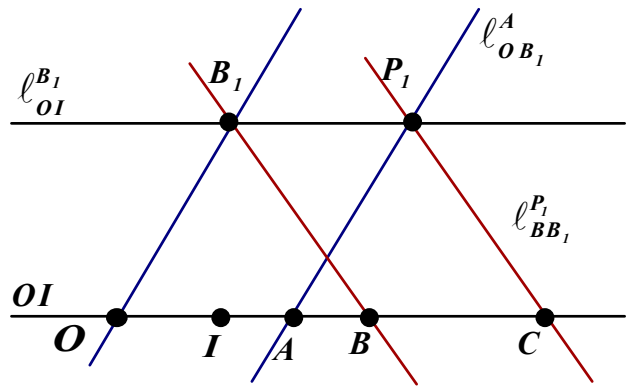

**Fig. 4 The additions configuration.**

The process of construct the points $C$, starting from two whatsoever points $A$, $B$ of the  line $OI$, is presented in the algorithm form

$$\left. \begin{array}{l} 1.\ B_1 \notin OI, \\ 2.\ \ell^{B_1}_{OI} \cap \ell^A_{OB_1} = P_1, \\ 3.\ \ell^{P_1}_{BB_1} \cap OI = C. \end{array} \right] \qquad (3)$$

In the process of construct the points $C$, besides pairs $(A, B)$ of points $A, B \in OI$ , is required and the selection of point $B_1 \notin OI$, which we call the *auxiliary point* to point $C$. This choice affects the position of point $C$ on the  line $OI$?

THEOREM 2.1. ***For every two points A, B $\in$ OI, algorithm*** (3) ***determines the a single point C $\in$ OI, which does not depend on the choice of hers auxiliary point $B_1$.***

*Proof*. Let it be $(A, B)$ a pair points of the  line $OI$.  According to (3), by selecting point $B_1$, construct the point $C$.

Now choose another point $B_2$. Then but according to (3), construct the analog point $C'$, in these conditions it takes view

$$\left.\begin{array}{l} 1.\ B_2 \notin OI, \\ 2.\ \ell_{OI}^{B_2} \cap \ell_{OB_2}^{A} = P_2, \\ 3.\ \ell_{BB_2}^{P_2} \cap OI = C'. \end{array}\right] \tag{3'}$$

We distinguish these four cases the position of points **A, B** in relation to fixed point **O** the fitting line **OI**.

**Case I. *A=B=O.*** In this case, by the choice of point $B_1$, according to (3) we have

$$P_1 = \ell_{OI}^{B_1} \cap \ell_{OB_1}^{O} = \ell_{OI}^{B_1} \cap OB_1 = B_1 \Rightarrow C = \ell_{OB_1}^{B_1} \cap OI = O;$$

whereas, from the choice of point $B_2$, according to (3') we have

$$P_2 = \ell_{OI}^{B_2} \cap \ell_{OB_2}^{O} = \ell_{OI}^{B_2} \cap OB_2 = B_2 \Rightarrow C' = \ell_{OB_2}^{B_2} \cap OI = O.$$

As a consequence (Fig. 5.*a*) we get

$$\boldsymbol{C=C'=O} \tag{4}$$

**Case II. *A=O≠B.*** In this case, by the choice of point $B_1$, according to (3) we have

$$P_1 = \ell_{OI}^{B_1} \cap \ell_{OB_1}^{O} = \ell_{OI}^{B_1} \cap OB_1 = B_1 \Rightarrow C = \ell_{BB_1}^{B_1} \cap OI = BB_1 \cap OI = B;$$

whereas, from the choice of point $B_2$, according to (3') we have

$$P_2 = \ell_{OI}^{B_2} \cap \ell_{OB_2}^{O} = \ell_{OI}^{B_2} \cap OB_2 = B_2 \Rightarrow C' = \ell_{BB_2}^{B_2} \cap OI = BB_2 \cap OI = B.$$

As a consequence (Fig. 5.*b*) we get

$$\boldsymbol{C=C'=B} \tag{5}$$

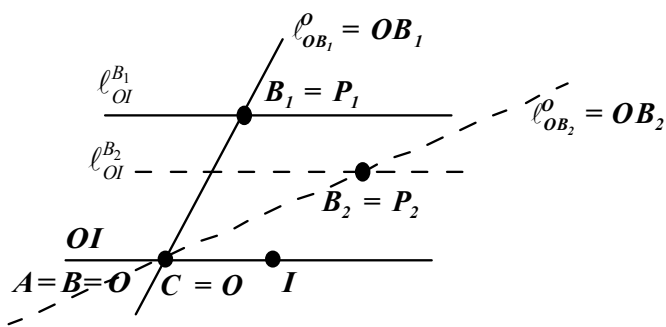

**Fig. 5.*a* independence of addition**

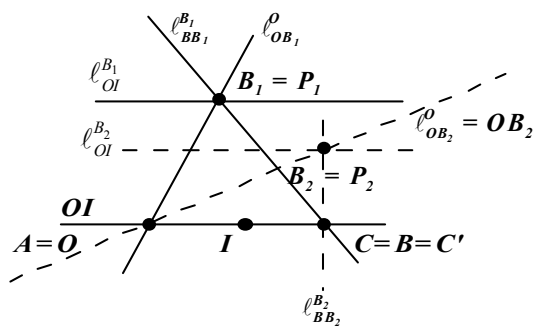

**Fig. 5.*b* independence of addition**

**Case III. *A≠O=B.*** In this case, by the choice of point $B_1$, according to (3) we have

$$P_1 = \ell_{OI}^{B_1} \cap \ell_{OB_1}^{A} \Rightarrow C = \ell_{OB_1}^{P_1} \cap OI = AP_1 \cap OI = A;$$

whereas, from the choice of point $B_2$, according to (3') we have

$$P_2 = \ell_{OI}^{B_2} \cap \ell_{OB_2}^{A} \Rightarrow C' = \ell_{OB_2}^{P_2} \cap OI = AP_2 \cap OI = A.$$

As a consequence (Fig. 5.*c*) we get

$$\boldsymbol{C=C'=A} \tag{5'}$$

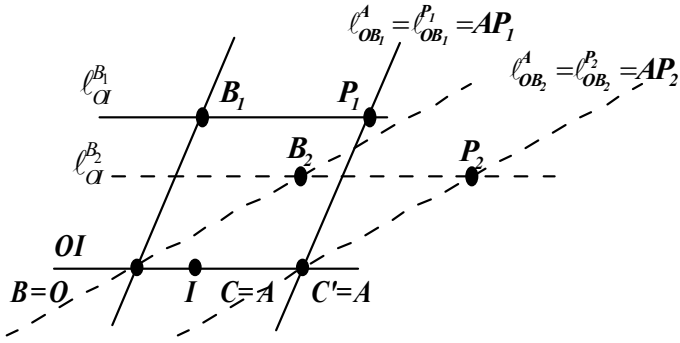 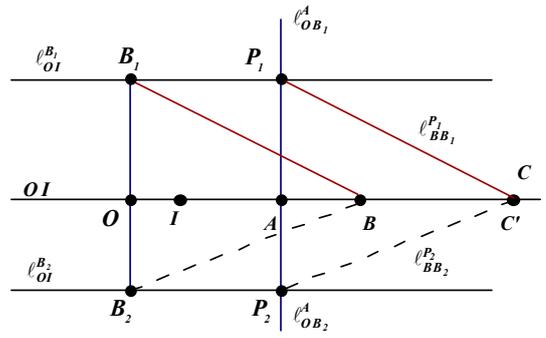

| **Fig. 5.**c **independence of addition** | **Fig. 5.**d **independence of addition** |

**Case IV**. $A \neq B \neq O$. Here we distinguish two sub-cases.

**a)** Points $O$, $B_1$, $B_2$ are collineary. In this case, by the choice of point $B_1$, according to (3) we have

$$P_1 = \ell_{OI}^{B_1} \cap \ell_{BB_1}^A \Rightarrow C = \ell_{BB_1}^{P_1} \cap OI;$$

whereas, from the choice of point $B_2$, according to (3') we have

$$P_2 = \ell_{OI}^{B_2} \cap \ell_{BB_2}^A \Rightarrow C' = \ell_{BB_2}^{P_2} \cap OI.$$

From (3) and (3') imply alsi, collinearity of the points $O$, $B_1$, $B_2$ imply collinarity of points $A$, $P_1$, $P_2$. Suppose now that $C \neq C'$ (Fig. 5.d).

We examine three-vertex $BB_1B_2$ and $CP_1P_2$. We note that $AP_1 = \ell_{OB_1}^A \parallel OB_1$, $P_2 \in AP_1$, $B_2 \in OB_1$, that imply $B_1B_2 \parallel P_1P_2$. But $C \in \ell_{BB_1}^{P_1} \parallel BB_1$, therefore $BB_1 \parallel CP_1$. From here, from axioms D1 of Desargues, results $B_2B \parallel P_2C$. On the other hand, $C' \in \ell_{BB_2}^{P_2}$, that imply $P_2C' \parallel B_2B$, which is parallel to $P_2C$. As a consequence $C' \in P_2C$. But $P_2C$ and $OI$ received in a single point, which imply $C=C'$, in contradiction with supposition that $C \neq C'$.

**b)** Points $O$, $B_1$, $B_2$ are non-collineary. In this case, by the choice of point $B_1$, according to (3) we have

$$P_1 = \ell_{OI}^{B_1} \cap \ell_{BB_1}^A \Rightarrow C = \ell_{BB_1}^{P_1} \cap OI;$$

whereas, from the choice of point $B_2$, according to (3') we have

$$P_2 = \ell_{OI}^{B_2} \cap \ell_{BB_2}^A \Rightarrow C' = \ell_{BB_2}^{P_2} \cap OI.$$

Suppose now that $C \neq C'$ (Fig. 5.e). From (3) and (3') we have, non-colinearity of points $O$, $B_1$, $B_2$ imply non-colinearity of the points $A$, $P_1$, $P_2$ .

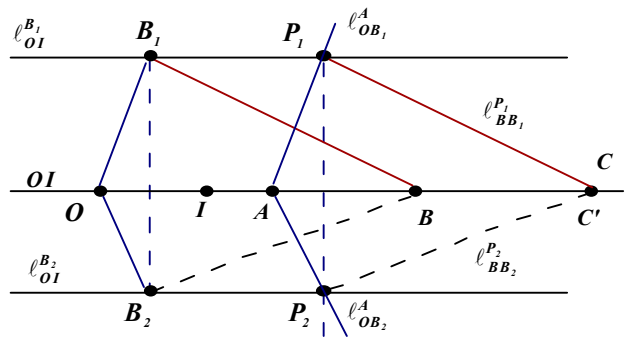

**Fig. 5.**e. **Independence of addition.**

We examine three-vertexes $OB_1B_2$ and $AP_1P_2$. We note that $AP_1 = \ell_{OB_1}^A \parallel OB_1$ and $AP_2 = \ell_{OB_2}^A \parallel OB_2$, therefore the axioms D1 results $B_1B_2 \parallel P_1P_2$. We examine now three-vertexes $BB_1B_2$ and $CP_1P_2$. The fact that $C \in \ell_{BB_1}^{P_1} \parallel BB_1$, imply $BB_1 \parallel CP_1$. From the above, we also $B_1B_2 \parallel P_1P_2$. Therefore, again by axioms D1 results $B_2B \parallel P_2C$. On the other hand, $C' \in \ell_{BB_2}^{P_2}$, that imply $P_2C' \parallel B_2B$. As a consequence $C' \in P_2C$, that imply $C=C'$, in contradiction with supposition that $C \neq C'$.

The above theorem creates the possibility of introduction of a binary operation, that we call the additions, in set of points to line $OI$, as follows.

Let us be $A$ and $B$ two whatsoever points of the line $OI$. I associate pairs $(A,B) \in OI \times OI$ point $C \in OI$, that determines algorithm (3). According to the preceding Theorems, point $C$ is determined in single mode by (3). Thus we obtain a application $OI \times OI \to OI$.

**Definition 2.1**. *In the above conditions, application*

$$+: OI \times OI \to OI,$$

*defined by* $(A, B) \quad C$ *for all* $(A, B) \in OI \times OI$ *we call the addition in* $OI$.

According to this Definitioni, can write

$$\forall A, B \in OI, \quad \begin{array}{l} 1.\ B_1 \notin OI, \\ 2.\ \ell_{OI}^{B_1} \cap \ell_{OB_1}^A = P_1, \\ 3.\ \ell_{BB_1}^{P_1} \cap OI = C. \end{array} \Bigg] \Leftrightarrow A + B = C. \tag{6}$$

## 3. GROUPOID ($OI$, +) IS COMMUTATIVE GROUP

With reference to cases I, II, III of Theorem 2.1, appears immediately true this
PROPOSITION 3.1. **Additions in $OI$ there areelement zero the point $O$**:

$$\forall\ A \in OI,\ O + A = A + O = A\ . \tag{7}$$

As well as worth and below propositions.

PROPOSITIONI 3.2. **Additions is commutative in $OI$**:

$$\forall\ A\ , B \in OI,\ A+B = B+A \tag{8}$$

*Proof.* In the case where $A=B$ the statement is evident, whereas when $A=O$ or $B=O$, propositions is tru goes according to (7). Stopped in case when $A, B \neq O$ and $A \neq B$. We mark $A+B=C$ and $B+A=C'$. Auxiliary point $B_1$ the sum $A+B$ and auxiliary point $A_1$ the sum $B+A$ we get the same (Fig.7). In this case, according to (6), we have

$$\begin{array}{l} 1.\ B_1 \notin OI, \\ 2.\ \ell_{OI}^{B_1} \cap \ell_{OB_1}^A = P_1, \\ 3.\ \ell_{BB_1}^{P_1} \cap OI = C. \end{array} \Bigg] \Leftrightarrow A + B = C \text{ and } \begin{array}{l} 1.\ A_1 \notin OI, \\ 2.\ \ell_{OI}^{A_1} \cap \ell_{OA_1}^B = P_2, \\ 3.\ \ell_{AA_1}^{P_2} \cap OI = C'. \end{array} \Bigg] \Leftrightarrow B + A = C'. \tag{6'}$$

It is clear that $A+B=B+A$ means that the points $C$ and $C'$ are the same points. For this use Proposition 1.2. Suppose now that $C \neq C'$ (Fig. 6). We examine trio collinary points $A$, $B$, $C$ and other trio of points collinary $B_1$, $P_1$, $P_2$, that are in parallel lines. According to (6'), $AP_1 \parallel BP_2$ and $BB_1 \parallel CP_1$. We are in conditions of little Pappus Theorems, thus resulting $CP_2 \parallel AB_1$, otherwise $CP_2 \parallel AA_1$. But from (6') have also $C'P_2 \parallel AA_1$, that imply $C=C'$, in contradiction with supposition that $C \neq C'$.

PROPOSITIONI 3.3. **Addition is associative in $OI$**:

$$\forall\ A, B, D \in OI,\ (A+B)+D = A+(B+D). \tag{9}$$

*Proof.* In the case where at least one of the point $A, B, D$ is $O$ proposition is tru according to (7), whereas when $A=D$, proposition is tru according to (8). Stopped in case the $A, B, D \neq O$ and $A \neq B \neq D$, (the reasoning is the same in other cases). Construct the first sum $(A+B)+D$. In this case (Fig. 7), according to (6), for $A+B$ have

$$1. \; B_1 \notin OI,$$
$$2. \; \ell_{OI}^{B_1} \cap \ell_{OB_1}^{A} = P_1, \quad \Bigg\} \Rightarrow A+B = \ell_{BB_1}^{P_1} \cap OI \;.$$
$$3. \; \ell_{BB_1}^{P_1} \cap OI = C.$$

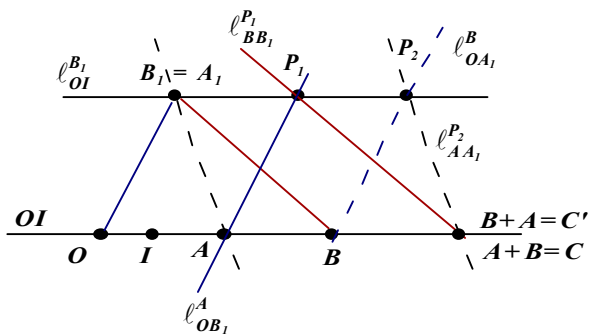

**Fig. 6.Commutative property**

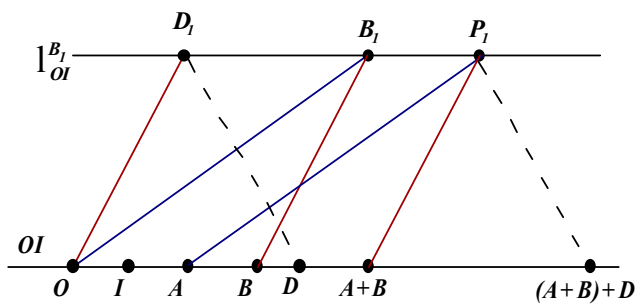

**Fig. 7. Associative property**

Construct line $\ell_{(A+B)P_1}^{O}$ and write down $D_1 = \ell_{(A+B)P_1}^{O} \cap OI.$ Select point $D_1$ as auxiliary points for construction of the sum $(A+B)+D.$ Then, according to (6) have

$$1. \; D_1 \notin OI,$$
$$2. \; \ell_{OI}^{D_1} \cap \ell_{OD_1}^{A+B} = P_1, \quad \Bigg\} \Rightarrow (A+B)+D = \ell_{DD_1}^{P_1} \cap OI \;. \qquad (*)$$
$$3. \; \ell_{DD_1}^{P_1} \cap OI = C.$$

On the order of same construct now sum $A+(B+D)$. In this case, we choose as auxiliary points for $B+D$ point $D_1$ (Fig. 8). With this, the role of point $P_1$ is the point $B_1$. By constructed line $\ell_{DD_1}^{B_1}$, according to (6), we find $B+D = \ell_{DD_1}^{B_1} \cap OI$. Whence imply that $(B+D)B_1 /\!/ DD_1$. Select now as auxiliary points for sum $A+(B+D)$ point $B_1$.

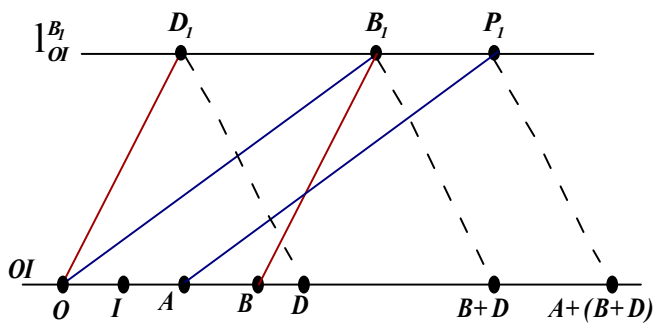

**Fig. 8. Associative property**

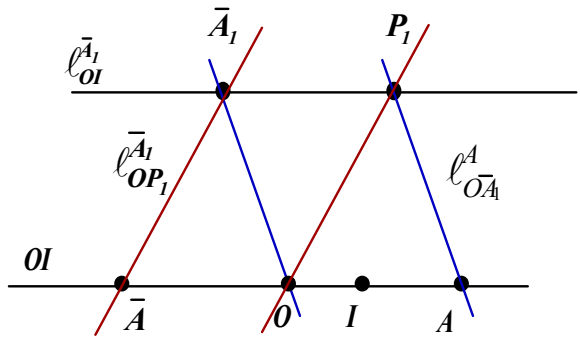

**Fig. 9.The inverse**

Then, according to (6) have

$$1. \; B_1 \notin OI,$$
$$2. \; \ell_{OI}^{B_1} \cap \ell_{OB_1}^{A} = P_1, \quad \Bigg\} \Rightarrow (A+B)+D = \ell_{(B+D)B_1}^{P_1} \cap OI \;.$$
$$3. \; \ell_{(B+D)B_1}^{P_1} \cap OI = C.$$

But $(B+D)B_1 /\!/ DD_1,$ that imply $\ell_{(B+D)B_1}^{P_1} = \ell_{DD_1}^{P_1}$. Eventually, according to (*), we have:

$$(A+B)+D = \ell_{DD_1}^{P_1} \cap OI = \ell_{(B+D)B_1}^{P_1} \cap OI = (A+B)+D.$$

PROPOSITION 3.4. **For every point in OI exists her right symmetrical according to addition**:

$$\forall \; A \in OI, \; \exists \overline{A} \in OI, \; A + \overline{A} = O. \tag{10}$$

*Proof.* We distinguish two cases: $A = O$ and $A \neq O$ .

If $A = O$ , then $\overline{A} = O$ , because, according to (7), $O + O = O$ .

If $A \neq O$ , requested points $\overline{A}$ such that

$$\left. \begin{array}{l} 1. \; \overline{A}_1 \notin OI, \\ 2. \; \ell_{OI}^{\overline{A}_1} \cap \ell_{O\overline{A}_1}^{A} = P_1, \\ 3. \; \ell_{\overline{A}\overline{A}_1}^{P_1} \cap OI = O. \end{array} \right]$$

Given this, we get initially a point $\overline{A}_1 \notin OI$ and construct line $\ell_{OI}^{\overline{A}_1}$ and then line $\ell_{O\overline{A}_1}^{A}$, which intersect at the point $P_1$ . Furthermore construct $OP_1$ and parallel with her by the points $\overline{A}_1$ line $\ell_{OP_1}^{\overline{A}_1}$. This last is not parallel with line $OI$, therefore awaits him at one point. It is clear that this point is the point of demanding $\overline{A}$, therefore $\overline{A} = \ell_{OP_1}^{\overline{A}_1} \cap OI$ (Fig. 9).

Propositions 3.1, 3.2, 3.3, 3.4 proved that is true this

THEOREM 3.1. **Groupoid (OI, +) is commutative(abeljan) Group**

*******